\documentclass[10pt,twocolumn,openbib]{article} 
\usepackage{graphicx}
\renewcommand{\baselinestretch}{1.0}
\newcommand{\Real}{\mbox{\boldmath $R$}}

\newcommand{\ed}{\end{document}}
\newcount\eqcount  
\def\eq{\global\advance\eqcount by1
        \eqno{\hbox{(\number\eqcount)}}}
\newcommand{\beqa}{\begin{eqnarray}}
\newcommand{\eeqa}{\end{eqnarray}}
\def\beqarr{\begin{eqnarray}}
\def\eeqarr{\end{eqnarray}}

\def\ux{\underline{x}}

\newcommand{\beq}{\begin{equation}}
\newcommand{\eeq}{\end{equation}}

\begin{document}
\newtheorem{theorem}{Theorem}
\newtheorem{proof}{Proof}
\renewcommand{\baselinestretch}{1.0}
%
%
\title{\bf Matching and Digital Control Implementation 
for Underactuated Systems}
\author{F. Andreev\thanks{Dept. of Mathematics, Visiting from Steklov Inst. of Math., St. Petersburg, Russia}, D. Auckly\thanks{Dept. of Mathematics}, L. Kapitanski$^{\ddagger}$, 
A. G. Kelkar$^{\S}$, and W. N. White\thanks{Department of Mechanical and Nuclear Engineering.\newline This work was partially supported 
by NSF Grant No. CMS-9813182.} \\Kansas State University\\Manhattan, KS 66506}
\date{}
\maketitle
 \thispagestyle{empty}
\abstract{This note describes two problems related to the digital 
implementation  of control laws in the infinite dimensional family of
matching control laws, namely state estimation and sampled data induced error.
The entire family of control laws is written for an inverted pendulum cart. 
Numerical simulations which include sampled data and a state estimator
are presented for one of the control laws in this family.}
\section{Introduction}
Several papers have been written recently regarding the control of 
nonlinear underactuated systems \cite{akw}-\cite{blm97}. 
Recall that an underactuated system, is a system 
with fewer control inputs than degrees of freedom. 
The idea that one should look for control laws such that the closed loop
system takes a particular form is common to all of these papers. The
particular form of the final equations is chosen so that there will be 
a natural candidate for a Lyapunov function. If the Lyapunov function 
attains a local  minimum at an isolated  point, then this point is 
a locally  asymptotically stable equilibrium of the continuous system. 

Digital implementation of continuous control laws introduces additional 
difficulties. It is not {\it a priori} clear that a method which produces good results in
the continuous case with full state feedback 
will continue to produce acceptable results
with state estimation and sampled data. For a digitally controlled system,
the data is collected and the control input is calculated at 
discrete moments of time. In addition, the full state cannot be directly measured 
and must be estimated based on observable data. 
Assume the continuous closed-loop system is modelled by 
\begin{equation}
\dot x\,=\, f(x,u(x)), \label{eq1}
\end{equation} 
where $u$ is a full state feedback control law. Let $\tau$ be the 
sample time,  let $x_k=x(\tau k)$, let $y_k=C(x_k)$ be the observerable 
data, and let $\ux_k$ be the estimated state at time $\,\tau k$. 
 A model of 
a  corresponding   digitally controlled system is 
\begin{equation}
  \begin{array}{c}
\dot x(t)\,=\, f(x(t),u(\ux(t)), \\
\ux(t)\,=\,\ux_k,\quad \hbox{for}\quad\tau k\le t<\tau (k+1),\\
\ux_{k+1}\,=\, g(y_k,\ux_k), 
   \end{array}
\label{e2}
\end{equation}
where $g$ is the state estimator. In practice, it is not system \ref{eq1}, but 
rather system (\ref{e2}), which must have an asymptotically stable equilibrium. 
There is a lower bound on the sampling time, $\tau$, dictated by the 
control apparatus. There is an upper bound on $\tau$ depending upon 
the state estimator and the control law.
 
For a simple linear system 
corresponding to $f(x,u)=Ax+Bu(x)$,  $y=Cx$, conditions to make (\ref{e2}) 
asymptotically stable are well known. For the continuous closed-loop system to  
be asymptotically stable,  there must exist  a matrix,
$G$ so that the eigenvalues of $A-GC$ all lie in the left half plane, 
and $\tau$ must be sufficiently small. 

For nonlinear systems, the situation is more complicated. 
One case, however, is easy to understand. Consider a nonlinear 
closed-loop system 
$
\dot x= f(x,u(x))
$
whose linearization, $\dot x= Ax +Bu(x)$, 
at the equilibrium $x=0$ is asymptotically stable such that  
there exists  a matrix
$G$ so that all eigenvalues of $A-GC$  lie in the left half plane. 
In this case, the continuous system
$$
\begin{array}{c}
\dot x= f(x,u({\ux})) \\
\dot{\ux} = f(\ux,u({\ux}))+GC(x-{\ux})
\end{array}
$$
will be locally asymptotically stable. 
In fact any system with correct linearization will do. 
For digital control, one may choose, 
for example, the following  system 
\begin{equation}
\begin{array}{c}
\dot x= f(x,u(\ux))\\
\ux(t)\,=\,\ux_k,\quad \hbox{for}\quad\tau k\le t<\tau (k+1),\\
\ux_{k+1}\,=\,A_d \ux_k + B_d u \ux_k + G_d C (x_k-\ux_k),
\end{array} \label{e3}
\end{equation}
where $A_d=\exp(\tau A)$, $B_d=\int_0^\tau \exp(-sA) ds B$  
and $G_d$ is chosen such that 
$\,spec\,(A_d-G_d C) = exp(\tau\,\,spec\,(A-GC))$. 
If the sampling time, $\tau$,  is sufficiently small, then the 
system (\ref{e3}) will be locally asyptotically stable.

\section{The Control Law}
In this section, we briefly recall a method for constructing an infinite
dimensional family of controls laws for many nonlinear systems. We will then use
the method to
derive a specific control law for an inverted pendulum cart system. The implementation
of this control law including a sample time and a state estimator will be discussed in the final section of this paper.

Let $\,Q\,$ denote a configuration space.
Let $g\in\Gamma(T^\ast Q\otimes T^\ast Q)$ be a metric.
Let $c,f:TQ\to TQ$ be fiber-preserving maps. We assume that $c(-X)=-c(X)$. 
Let $V:Q\to{\Real}$. The
system that we consider is
\[
 \nabla_{\dot\gamma}\dot\gamma 
+ c(\dot\gamma) +\ grad_\gamma V
   =f(\dot\gamma).   
\]
Let $P\in\Gamma(T^\ast Q\otimes TQ)$ be a $g$-orthogonal projection. 
We consider the situation where a constraint $P(f) = 0$ is imposed.
A system is called underactuated if $P \neq 0$.
In order to describe the final control law, we will use several other
variables. The variable $\widehat g\in\Gamma(T^\ast Q\otimes T^\ast Q)$ will
be a metric,
 $\widehat c:TQ\to TQ$ will be a fiber-preserving map, $\widehat V$
will be a real-valued function, and $\lambda \in\Gamma(T^\ast Q\otimes TQ)$
will be a $g$-self adjoint map. One first solves the equations 
$$ 
\nabla g\lambda\big|_{\hbox{Im}\ P^{\otimes 2}}=0,
   $$
for $\lambda |_{\hbox{Im}\ P}$. Then one solves 
$$
L_{{}_{\lambda PX}}\widehat g = L_{{}_{PX}} g 
$$
(this is a slight rewrite of equation (1.12) of our previous paper \cite{akw}),
$$
L_{{}_{\lambda PX}}\widehat V = L_{{}_{PX}} V 
$$
(this is equation (1.13) of our previous paper \cite{akw}),
then after solving,
$$ 
P(c(X)-\widehat c(X))=0,  
$$
the control input will be given by:

\begin{equation}
f(X)\equiv \nabla_X X-\widehat \nabla_X X+\ grad_\gamma V
   -\ \widehat {grad}_\gamma \widehat V +c(X)-\widehat c(X)
\label{e4}
\end{equation}

We now apply the above method to the inverted pendulum cart depicted in 
Figure 1.

\begin{figure}[htbp]
\centerline{\includegraphics[width=2.5in]{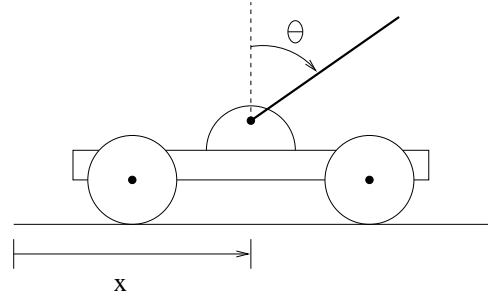}}
\caption{Inverted pendulum cart}
\end{figure}
With appropriate scaling, the metric $g$ is given by 
$g=d\theta^2 +2 b \cos(\theta)\,dx\,d\theta + dx^2$, where $b$ is a physical parameter, $0<b<1$. The potential energy is given by $V=\cos(\theta)$. 
Since no torques can be applied directly to the pendulum, 
$P=(b\cos(\theta)\,dx+d\theta)\otimes {\partial/\partial\theta}$ is 
the orthogonal projection onto the direction ${\partial/\partial\theta}$. 
Assuming that there is no dissipation, $c=0$.

Let $\theta$ be the coordinate with index $1$, and $x$ be the coordinate with
index $2$. Writing $\lambda PX=\sigma {\partial\over\partial \theta}
+\mu {\partial\over\partial x}$, where $\sigma$ and 
$\mu$ are yet to be found, the $\lambda$-equation may be rewritten as
$$
{\partial\over\partial \theta} (\sigma+b \cos(\theta) \mu)+2b\sin(\theta)\,
\mu=0,
$$
$$
{\partial\over\partial x} (\sigma+b\cos(\theta)\mu)=0. 
$$
For these equations to be consistent 
the following compatibility condition must hold:
$$
{\partial\over\partial x} (\sin(\theta)\mu)=0.
$$
This implies that $\mu$ is a function of $\theta$. The second 
$\lambda$-equation implies that $\sigma$ is a function of $\theta$.
The first $\lambda$-equation then becomes an ODE which may be solved for
$\sigma$ giving,
$$
\sigma(\theta)=\sigma_0+b\mu_0-b\cos(\theta)-2b\int^{\theta}_0 
\sin(t) \mu(t) \, dt.
$$
Before solving the $\widehat g$-equation, it is helpful to solve,
$$
\sigma {\partial y\over\partial \theta}+\mu {\partial y\over\partial x} =0.
$$
Using the method of characteristics, we find,
$$
y=x-\int^{\theta}_0 {\mu(t)\over\sigma(t)} dt.
$$
The $\widehat g$-equation may be rewritten as
$$
\sigma {\partial\widehat g_{11}\over\partial \theta}+
\mu {\partial\widehat g_{11}\over\partial x}
+2({\partial\sigma\over\partial \theta}
-{\sigma\over\mu}{\partial\mu\over\partial \theta})\widehat g_{11}
+2{{\partial\mu\over\partial \theta}\over\mu}=0.
$$
Let $\bar\sigma$ and $\bar\mu$ 
be $\sigma$ and $\mu$ considered as functions of $\theta$ and $y$, i.e.,
$$
\bar\sigma( \theta,y( \theta,x))=\sigma(\theta,x),\qquad
\bar\mu(\theta,y(\theta,x))=\mu(\theta,x). 
$$
The solution to the $\widehat g$-equation is then given explicitly by
$$
\widehat g_{11}( \theta,x)={\mu^2\over\sigma^2}\left[
-2\,\int_0^{\theta}{\bar\sigma\over\bar\mu^3}
{\partial\bar\mu\over\partial \theta}\,d \theta\,\big|_y
+h(y)\right], 
$$
where $h(y)$ is an arbitrary function of a single variable. Using the 
definition of $\lambda$, we have
$$
\widehat g_{12} ={1\over\mu} (1-\sigma \widehat g_{11}),
$$
$$
\widehat g_{22} ={1\over\mu} (b\cos(\theta)-\sigma \widehat g_{12}). 
$$
Using integration by parts and the first $\lambda$-equation, we can simplify
the integral appearing in $\widehat g_{11}$. Explicitly,
\[
\widehat g_{11}( \theta,x)={1\over\sigma}-{\sigma_0\mu^2\over\mu_0^2\sigma^2}
-{b\cos(\theta)\mu\over\sigma^2}+{b\mu^2\over\mu_0\sigma_0^2}
+{\mu^2\over\sigma^2}h(y).
\]

The function $\widehat V$ satisfies the equation
$$
\sigma {\partial \widehat V\over\partial \theta}+
\mu {\partial \widehat V\over\partial x} =-\sin(\theta).
$$
Considering $\widehat V$ as a function of $\theta$ and $y$, this becomes an ODE.
The resulting expression for $\widehat V$ is:
$$
\widehat V(\theta, x)=w(y(\theta,x))-\int^\theta_0 {\sin(t)\over\sigma(t)} \, 
dt,
$$
here $w(y)$ is an arbitrary function.
Finally, the solution to the $\widehat c$-equation is given by:
$$
\widehat c^1=-b\cos(\theta)(\widehat c^2 -c^2)+c^1,
$$
where $\widehat c^2$ is an arbitrary function which is odd in the velocities. 
The final control law is given by equation (\ref{e4}).

For the inverted pendulum in our lab, the parameter, $b$, is .238, and $c=0$.
(In practice, there is some dissipation in the base of the cart, but this may be
directly counteracted by a term in the control law. The dissipation in the
joint holding the pendulum is really negligible.) The value of $b$ is  
different from the value that was used in the simulations in \cite{akw}. This is
because we are
now taking into account additional contributions to the mass of the base
of the cart and to the mass of the pendulum. 

The control law studied in \cite{akw} stabilized a wide range of initial conditions, 
however, it was found to 
be underdamped for small initial conditions. This was not dissapointing
because the arbitrary functions in our control law were chosen for algebraic
simplicity, and not for specific engineering goals. For this paper, we
decided to use step functions in place of some of the constants used previously.
By using step functions we hoped to reduce the number of parameters to 
something which would be reasonable to analyze. Our plan was to try to blend
the nonlinear control law which worked well for large disturbances with
one which would linearize to the linear control law that worked well for
small initial conditions. In particular, we took
$$
\mu(\theta)= \left\{ \begin{array}{ll} \mu_0 \cos(\theta), & \hbox{for $|\theta|\leq\theta_L$}\\
\mu_{\infty} \cos(\theta), & \hbox{otherwise}
\end{array} \right. .
$$
$$
h(y)= \left\{ \begin{array}{ll} h_0, & \hbox{for $y\leq y_L$}\\
h_{\infty} , & \hbox{otherwise}
\end{array} \right. 
$$
$$
w(y)= \left\{ \begin{array}{ll} \frac12 w_0 y^2, & \hbox{for $y\leq y_L$}\\
\frac12 w_{\infty} y^2, & \hbox{otherwise}
\end{array} \right. 
$$
The entire motivation for this method is that
$
\widehat H(\dot \gamma)=\frac12 \widehat g(\dot\gamma, \dot\gamma)+\widehat V(\gamma)
$, 
is a natural candidate for a Lyapunov function for the closed loop system.
The time derivative of $\widehat H$ is, $-\widehat g(\widehat c(X),X)=
(\det \widehat g)\cdot \widehat c^2\cdot (\mu_0\cos\theta\dot\theta-\sigma_0\dot x)$. Thus taking $\widehat c^2=\Phi(\mu_0\cos\theta\dot\theta-\sigma_0\dot x)$
will insure that $\widehat H$ is never increasing. We take, 
$$
\Phi(\theta)= \left\{ \begin{array}{ll} \Phi_0 , & \hbox{for $|\theta|\leq\theta_L$}\\
\Phi_{\infty} , & \hbox{otherwise}
\end{array} \right. .
$$
With these choices, the function $\sigma$ will take the form,
$$
\sigma(\theta)= \left\{ \begin{array}{ll} \sigma_0 , & \hbox{for $|\theta|\leq\theta_L$}\\
\sigma_{\infty} , & \hbox{otherwise}
\end{array} \right. ,
$$
where $\sigma_{\infty}=\sigma_0+b(\mu_0-\mu_{\infty})\cos^2(\theta_L)$.
We guessed $\theta_L=.3$ and $y_L=15$. Values which stabilize a large
region are: 

$\sigma_{\infty}=-.05,~~ \mu_{\infty}=9.9,~~ w_{\infty}=1.5,$
 $\Phi_{\infty}=.75,$ and $h_{\infty}=.03.$

These are a slight modification of the values given in our previous paper,
since we are using a slightly diferent value of $b$.
To compute the appropriate values for the remaining constants, we write
the linear control input as: $u_l=g(f,{\partial\over\partial x})=
p_1\theta+p_2x+d_1\dot\theta+d_2\dot x$. Setting ${\partial u\over\partial \theta}|_0=p_1$, ${\partial u\over\partial x}|_0=p_2$, 
$d_1\sigma_0+d_2\mu_0=0$, and $\sigma_{\infty}=\sigma_0+b(\mu_0-\mu_{\infty})\cos^2(\theta_L)$ determines
four of the remaining parameters. The final parameter is determined 
by the condition, ${\partial u\over\partial \dot\theta}|_0=d_1$. 
The resulting parameters are: 

$\sigma_0=-1.59$, $\mu_0=17$, $w_0=.00296$, 
$\Phi_0=1.48$, and $h_0=.0081.$

Numerical results comparing the control law described above with the
linear control law are presented in Figs. 2 through 5 below. The small initial 
conditions were $\theta_0=0.4$, $x_0=0$, $\dot\theta_0=0$, and $\dot x_0=0$. The large 
initial conditions were $\theta_0=1.1$, $x_0=0$, $\dot\theta_0=0$, and $\dot x_0=0$.

\begin{figure}[htbp]
\centerline{\includegraphics[width=3.0in]{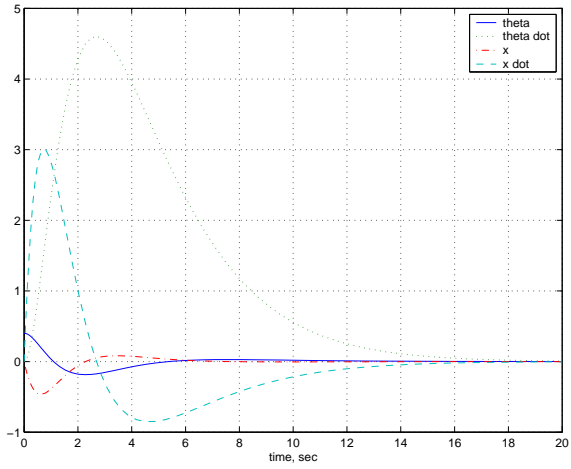}}
\caption{Linear control law with small initial conditions}
\end{figure}

\begin{figure}
\centerline{\includegraphics[width=3.0in]{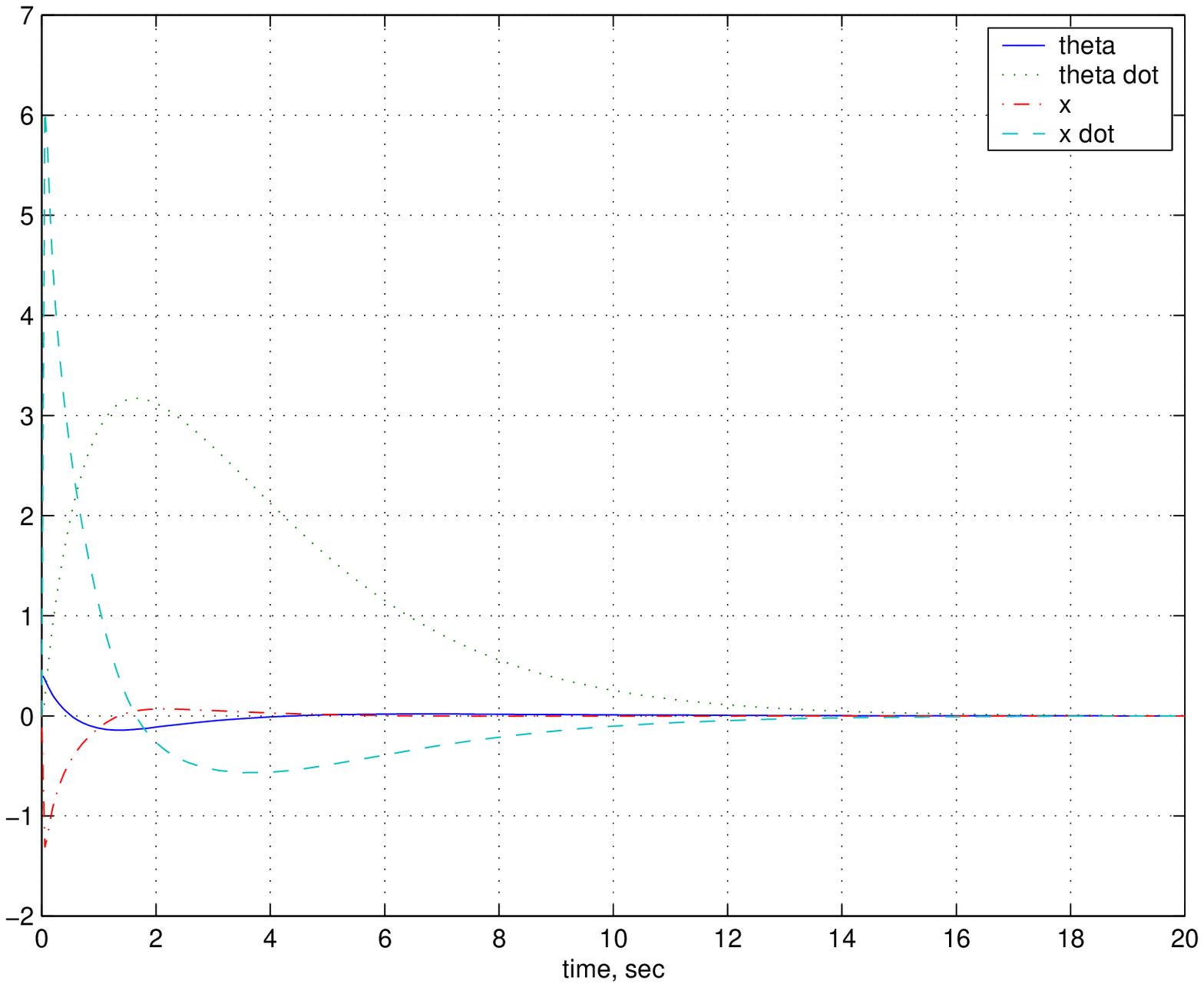}}
\caption{Nonlinear control law with small initial conditions}
\end{figure}

\begin{figure}
\centerline{\includegraphics[width=3.0in]{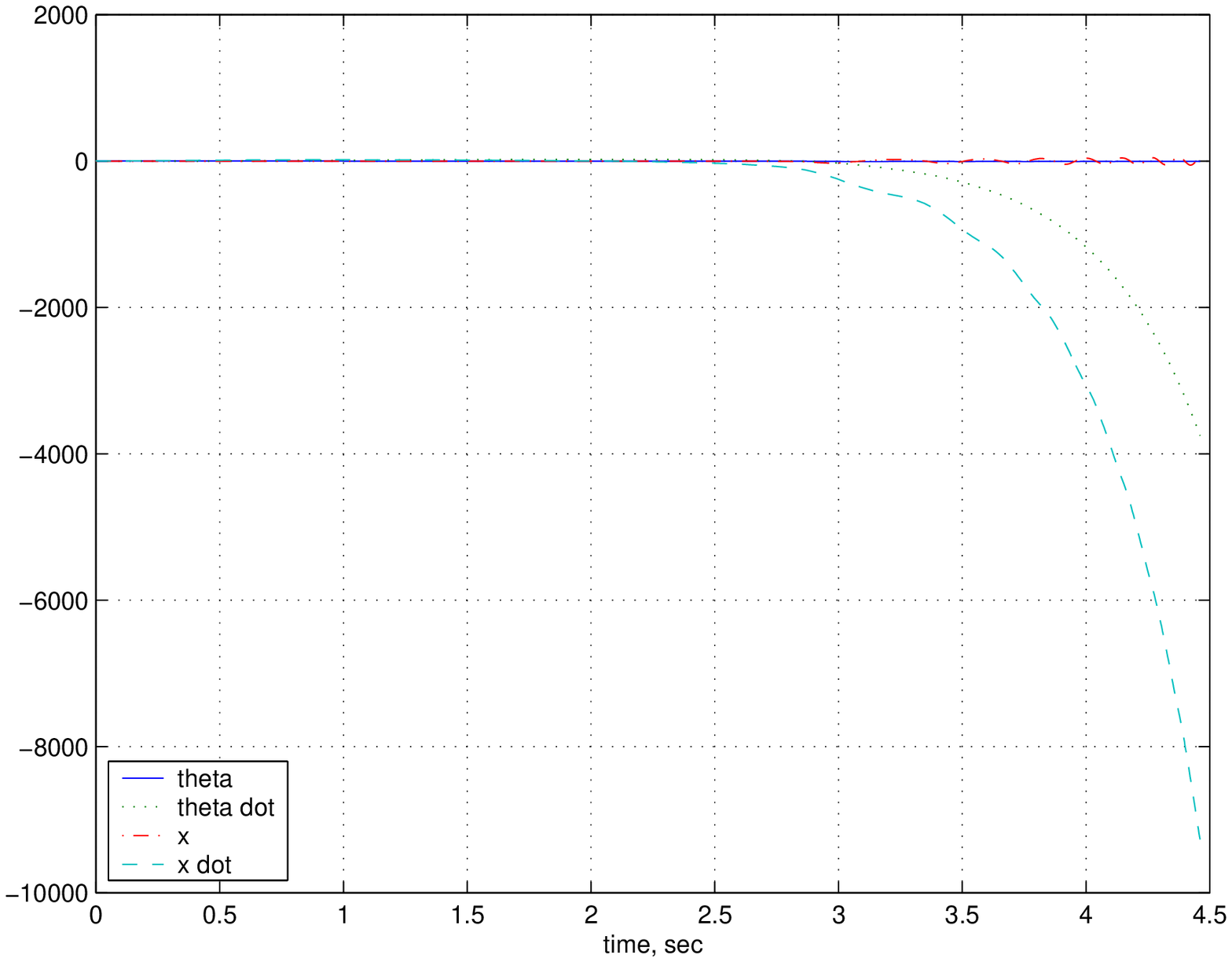}}
\caption{Linear control law with large initial conditions}
\end{figure}

\begin{figure}
\centerline{\includegraphics[width=3.0in]{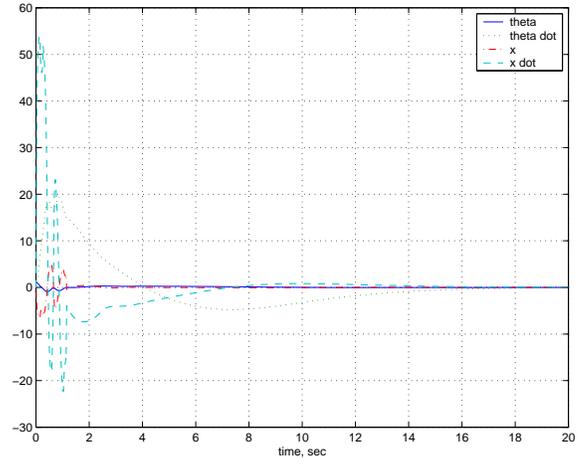}}
\caption{Nonlinear control law with large initial conditions}
\end{figure}

\section{Implementation}
The inverted pendulum cart in our lab cannot directly observe the
velocity or angular velocity of the cart. Thus, full state feedback is
not possible and the modifications needed to implement the control
law based upon the output of a non-trivial
linear observer, $C$, must be considered.
For our first test of a digital control system implementing the control law
described above, we took an estimator of the form (\ref{e3}) with $\tau=0.0143$,
\[
A_d = \left( \begin{array}{llll} 1&0&0.0143 & 0 \\
                         0&1&0&0.0143 \\
                         0.0151&0&1&0 \\
                         -.0036&0 &0&1 \end{array} \right),
\]
\beqa
 B_d &=& \left( \begin{array}{l} 0 \\
                         0 \\
                         -0.0036 \\
                         0.0151 \end{array} \right), 
C = \left( \begin{array}{llll} 1&0&0 & 0 \\
                         0&1&0&0 \\
                         \end{array} \right), \nonumber\\
G_d &=& \left( \begin{array}{ll} 0.168 & 0 \\
                         -.0001&0.165 \\
                         0.509&0 \\
                         -.0039&0.473 \end{array} \right), 
x=\left( \begin{array}{l} \theta \\
                         x \\
                         \dot\theta \\
                         \dot x \end{array} \right) \nonumber
\eeqa
Results from our numerical computations of this larger system are displayed 
in Figures 6 and 7. The initial estimated state was chosen to be the same as the initial conditions. System
 (\ref{e3}) was unstable with the large initial conditions for both the linear and nonlinear control law. 


\begin{figure}[htbp]
\centerline{\includegraphics[width=3.0in]{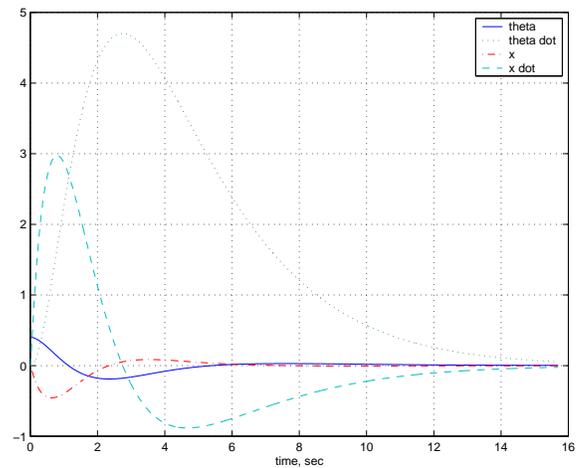}}
\caption{Linear control law with small initial conditions
state estimation and sampled data}
\end{figure}

\begin{figure}[htbp]
\centerline{\includegraphics[width=3.0in]{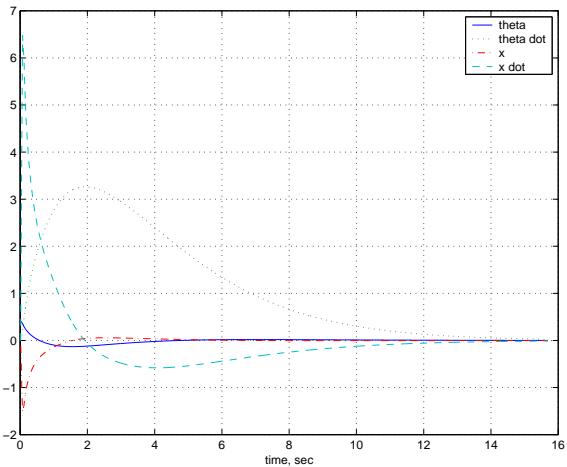}}
\caption{Nonlinear control law with small initial conditions, 
state estimation and sampled data}
\end{figure}
\section{Conclusions}
Looking for control laws so that the closed loop system takes some specified form appears 
to be a promising idea in nonlinear control theory. There are 
however many issues which have not been fully resolved. One must first decide
what it means to say that one control law is better than another. With only
an intuitive idea of what is "better" we would argue that a control law derived via
the matching equations works ``better" for an inverted pendulum cart than a
linear control law. This brings up the question of finding control laws of
this type. Such control laws are usually described as solutions to a system
of partial differential equations. Just guessing a solution based on the
form of the equations is not a very satisfactory solution to the problem.
The general solution to the matching equations may be found for systems
with two degrees of freedom. If there is some symmetry present, it is also
possible to find solutions to the matching equations. This leaves open the
problem of finding such control laws for systems with more degrees of freedom
in the absence of symmetry. The general matching equations have many solutions.
This means that one must have some method for picking a good solution to 
the matching equations.  

Assuming that all of these questions have been answered, one must still
come up with satisfactory answers to the main questions discussed in this
paper: what state estimator should be used, and why will the closed loop
system be stable when only sampled data is used. Mathematically, it is
well known that the resulting system will be locally stable if the 
linearization about the equilibrium is stable and the sample time is
sufficiently small. Perhaps some numerical and experimental tests will
shed some light on the correct choice for a state estimation scheme. Given
the promising results of the matching control law applied to the inverted
pendulum cart, and the wide array of open questions, this is a fertile area
for future research. 



\end{document}